\documentclass[twoside, 10pt]{article}
\usepackage{mathrsfs}

\usepackage{mathrsfs,amsfonts,amsmath}
\RequirePackage{ifthen,calc}
\usepackage{theorem}
\usepackage[dvips]{color}

 \renewcommand{\baselinestretch}{1.0}

 \catcode`@=11
 \def\@evenhead{\hbox to\textwidth{\footnotesize\rm\thepage \hfill
  {\it Yuqiang LI}}} 

 \def\@oddhead{\hbox to \textwidth{\footnotesize{\it
   Fluctuation limits of strongly degenerate branching systems } \hfill\thepage}}

 \renewcommand{\section}{\makeatletter
 \renewcommand{\@seccntformat}[1]{{\csname the##1\endcsname.}\hspace{0.45em}}
 \makeatother \@startsection
{section}
{1}
{0pt}
{\baselineskip}
{0.5\baselineskip}
{\normalsize\bfseries\mathversion{bold}}}

\catcode`@=12

\newtheorem{thm}{\noindent Theorem}[section]
\newtheorem{lem}{\noindent Lemma}[section]

\theoremheaderfont{\normalfont\bfseries} \theorembodyfont{\slshape}
{\theorembodyfont{\rmfamily}} {\theorembodyfont{\rmfamily}
}
{\theorembodyfont{\rmfamily} } {\theorembodyfont{\rmfamily}
}
{\theorembodyfont{\rmfamily} } {\theorembodyfont{\rmfamily}
\newtheorem{rem}{\noindent Remark}[section]}
{\theorembodyfont{\rmfamily} } {\theorembodyfont{\rmfamily} } {\theorembodyfont{\rmfamily}
}

 \def\beqlb{\begin{eqnarray}}\def\eeqlb{\end{eqnarray}}
 \def\beqnn{\begin{eqnarray*}}\def\eeqnn{\end{eqnarray*}}
 \numberwithin{equation}{section}
\def\qed{\hfill$\square$\smallskip}
\def\R{{\mathbb R}}
\def\bfE{{\mathbb{E}}}

\def\e{\mathrm{e}}

\def\1int{\int_{0+}^\infty}
\def\0int{\int_0^\infty}
\def\d{\mathrm{d}}

\begin{document}

\title{\LARGE\bf Fluctuation limits of  strongly degenerate branching systems
}
\author{ Yuqiang Li \renewcommand{\thefootnote}{\fnsymbol{footnote}}\footnotemark[1]
\\ \small School of Finance and Statistics, East China Normal University,
\\ \small Shanghai 200241, P. R. China.
}
\date{}

\maketitle
\renewcommand{\thefootnote}{\fnsymbol{footnote}}\footnotetext[1]
{Research supported partly by NSFC grant (10901054).}

\centerline{\textbf{Abstract}}

\noindent{Functional limit theorems for scaled occupation time fluctuations
of a sequence of
generalized branching particle systems in $\R^d$ with anisotropic
space motions and strongly degenerated splitting abilities are studied in the cases of critical and intermediate dimensions. The
results show that the limit processes are time-independent measure-valued Wiener processes with simple spatial structure.}

\smallskip


\smallskip


\medskip

\noindent\textbf{Running Title}: Fluctuation Limits of strongly degenerate
 branching systems

\renewcommand{\baselinestretch}{1.2}

\normalsize

\bigskip

\section{Introduction}
Consider a kind of generalized branching particle systems in $\R^d$.  Particles
start off at time $t=0$ from a Poisson random field with Lebesgue
intensity measure $\lambda$ and evolve independently. They move in $\R^d$ according to a L\'evy process
$$
\vec{\xi}=\{\vec{\xi}(t),
 t\geq0\}=\{(\xi_1(t),\xi_2(t),\cdots,\xi_d(t)), t\geq 0\}
 $$
with independent stable components as in \cite{PT69}, i.e. for
every $0<k\leq d$, $\xi_k=\{\xi_k(t), t\geq0\}$ being a symmetric
$\alpha_k$-stable L\'{e}vy process and $\xi_1,\cdots,\xi_d$
independent of each other. In addition, the
particles split at a rate $\gamma$ and the branching law at age $t$ has the generating function
 $$g(s,t)=\Big(1-\frac{\e^{-\delta t}}{2}\Big)
 +\e^{-\delta t}\frac{s^2}{2}, \qquad 0\leq s\leq1,\  t\geq 0.$$
Intuitively, in this model, the particles' movement in different
direction is controlled by different mechanism and their probability of
splitting new particles declines with the rate $\delta$ as their ages increase.
It is easy to see that when $\delta=0$, this model is similar to a classical $(d,\alpha,\beta)$-branching particle system with $\beta=1$ except that the moving mechanism is the anisotropic stable L\'{e}vy processes $\vec{\xi}$ rather than a  symmetric $\alpha$-stable L\'{e}vy process. Li and Xiao \cite{LX10} called this model as a $(d, \vec{\alpha},
\delta,\gamma)$-degenerate branching particle system, where $\vec{\alpha}:=(\alpha_1,\cdots,\alpha_d)$.  Let
 $\bar{\alpha}:=\sum_{k=1}^d1/\alpha_k.$
When $\bar{\alpha}>2$, $\bar{\alpha}=2$ and
$\bar{\alpha}\in(1,2)$, the corresponding dimension of the space
is referred to as the large dimension, critical dimension and
intermediate dimension, respectively.

Motivated by the work on occupation time fluctuations of classical branching particle systems and the work on construction of anisotropic random fields (see, for example, \cite{BMS07}), Li and Xiao \cite{LX10} explicitly studied the functional limits of occupation time fluctuations of the models. Observe that a fixed $(d, \vec{\alpha}, \delta,\gamma)$-branching particle system with $\delta>0$ will go to local extinction as time elapses because of the sub-critical branching laws at positive ages. They \cite{LX10} borrowed the idea of nearly critical branching processes (see \cite{IPZ05,Li09,ST94}) and considered a sequence of $(d,\vec{\alpha}, \delta_n,\gamma)$-models with $\delta_n\to0$ as $n\to\infty$. More precisely, let $N_n(s)$ denote the empirical measure of the $(d, \vec{\alpha},\delta_n,\gamma)$-degenerate branching particle system at
time $s$, i.e. $N_n(s)(A)$ is the number of particles in the set
$A\subset\R^d$ at time s.  They studied the limit of a sequence of scaled occupation time fluctuations,
 \beqlb\label{s1-1}
 X_n(t)=\frac{1}{F_n}\int_0^{nt}(N_n(s)-f_n(s)\lambda)\d s,
 \eeqlb
where $F_n$ is a scaling constant and
  \beqlb\label{s1-2}
  f_n(s):=\bar{f}_n(s)\e^{-\delta_n s}:=\Big[1+\frac{\delta_n}{\gamma-\delta_n}(1-\e^{-(\gamma-\delta_n)s})\Big]\e^{-\delta_n s},
  \eeqlb
under the assumption $n\delta_n\to\theta\in[0,\infty)$ which is referred to weak degeneration, and proved that in the cases of critical and intermediate dimensions the limit processes have complicated temporal structures and in the case of large dimensions, the limit processes own simple temporal but anisotropic spatial structures.


The purpose of this paper is to continue the discussion of functional limits of (\ref{s1-1}) under the assumption that $n^\kappa\delta_n\to\theta\in(0,\infty)$ for some $\kappa\in(0, 1)$, which is referred to strong degeneration. We  focus on the cases of critical and intermediate dimensions in this paper. The main methods used in this situation is same as that in Li and Xiao \cite{LX10}, which was formulated and developed by Bojdecki {\it et al} in their serial papers (\cite{BGT061}-\cite{BGT072}), except some complexities and
differences from the strong degeneration.  We find that the limit processes in any
positive time interval are time-independent measure-valued Wiener processes, which always have the form $C\lambda\xi$, where $\xi$ is a standard normal random variable, $\lambda$ is the Lebesgue measure in $\R^d$ and $C$ is a non-random constant. By comparison with the corresponding results in Li and Xiao \cite{LX10}, the current limit processes are simpler (please see Remark 2.1 in Section 2 for more details). To save the space of this paper, we leave the study on the case of large dimensions elsewhere because the potential limit processes deserve further investigations. In addition, we remark that there are few results under the assumption that $\delta_n\to 0$ and $n^\kappa\delta_n\to\infty$ for any $\kappa>0$.

There is much literature related to the field of fluctuations of branching particle systems. Iscoe \cite{I86} studied the
single time limit theorem of occupation time of the  $(d,
\alpha,\beta)$- superprocess which in essence is the limit process
of the classical $(d, \alpha,\beta)$-branching system, and got different
limits depending on the relations between $d, \alpha, \beta$. Hong [16] also considered the superprocess case and proved the
convergence of finite-dimensional distributions of real processes without the tightness for a fixed test
function. Recently Bojdecki et al in their series of papers, such as \cite{BGT061}-\cite{BGT082}, studied the functional limits of occupation time fluctuations of a fixed classical  $(d, \alpha,\beta)$-branching system which is different from the setting in Li and Xiao \cite{LX10} and this paper. For more literature we refer to \cite{DP99,DG01,L98} and the references therein.

Without other statement, in this paper, we use $K$ to denote an
unspecified positive finite constant which may not necessarily be
the same in each occurrence. In addition, since this paper and
\cite{LX10} discuss the same branching systems in different assumptions,
in order to shorten the length of this paper, we will omit some common
inferences and calculations and refer to \cite{LX10}.

The remainder of this paper is organized as follows. Section 2 contains the
main results of this paper and some auxiliary results and
formulas used in the proofs of the main results. In Section 3 we prove the main results.

\section{Main results}

Consider a sequence of $(d,
\vec{\alpha},\delta_n,\gamma)$-degenerate branching particle system. The particles' spatial movement is described by $\vec{\xi}_n$. We assume
that $\{\vec{\xi}_n, n \ge 1\}$ is a sequence of identically
distributed $\R^d$-valued L\'evy processes with $\alpha_k$-stable
components ($1\le k \le d$). The distribution of $\vec{\xi}_n$ is
completely determined by its characteristic function
\begin{equation}\label{Eq:ChF}
{\mathbb E}\Big(\e^{i \langle z, \vec{\xi}_n(t)\rangle}\Big) =
\e^{- t \sum_{k=1}^d |z_k|^{\alpha_k}}, \quad z \in \R^d.
\end{equation}
Obviously, for any $n>0$, $\vec{\xi}_n$ is a
time-homogeneous Markov process on $\R^d$. Since $\vec{\xi}_n$ has
the same distribution for all $n$, we denote its semigroup by
$\{T_t\}_{t\geq0}$, i.e.,
 $$T_sf(x):=\bfE(f(\vec{\xi}_n(t+s))|\vec{\xi}_n(t)=x),$$
for all $s, t\geq0$, $x\in\R^d$ and bounded measurable functions
$f$ on $\R^d$.  In order to avoid misunderstanding, in case of necessity we write $T_sf(x)$ by $T_s(f(\cdot))(x)$.

In this paper, we always let $N_n(t)$ be the empirical measure of the $(d,
\vec{\alpha},\delta_n,\gamma)$-model for every $n\geq 1$. From Section 2 in Li and Xiao \cite{LX10}, we know that the scaled occupation
time fluctuations of $(d,
\vec{\alpha},\delta_n,\gamma)$-models are defined as follows.
 \beqlb\label{s2-1}
 \langle X_n(t),\phi\rangle=\frac{1}{F_n}\int_0^{nt}\langle N_n(s)-f_n(s)\lambda,
 \phi\rangle\d s,
 \eeqlb
for every $\phi\in\mathcal{S}(\R^d)$, the space of smooth rapidly decreasing
functions, where $F_n$ is a suitable
scaling parameter, $f_n(s)$ same as (\ref{s1-2}) and $\langle\mu, f\rangle=\int f\d\mu$ for any measure $\mu$ and
any integrable function $f$ on $\mu$.
Furthermore,
  \beqlb\label{s2-2}
  \bfE(\langle N_n(s), \phi\rangle|N_0=\epsilon_x)=f_n(s)T_s\phi(x),
  \eeqlb
for all $n\geq0, x\in\R^d$ and $\phi\in\mathcal{S}(\R^d)$. Here
$\epsilon_x$ denotes the unit measure concentrated at $x\in\R^d$.

Below, we assume  that there is a constant $\theta\in(0,\infty)$ such that $n^{\kappa}\delta_n\to\theta$ for some $k\in(0, 1)$ as $n\to\infty$.
Let $\widehat{\phi}(z)$ $(z\in\R^d)$ be the Fourier
transform of function $\phi\in L(\R^d)$, i.e.,
 $\widehat{\phi}(z)=\int_{\R^d}\e^{i\langle x, z\rangle}\phi(x)\d x$, and $\mathcal{S}'(\R^d)$ the dual
space of $\mathcal{S}(\R^d)$. Recall that $\vec{\alpha}:=(\alpha_1,\cdots,\alpha_d)$ and $\bar{\alpha}:=\sum_{k=1}^d1/\alpha_k.$
The main results of this paper read as follows.

\begin{thm}\label{s2-thm-1}
When $\bar{\alpha}=2$, let $F_n^2=n^{\kappa}\ln
n$. Then for any $\varepsilon>0$, $X_n\Rightarrow C_1\lambda\zeta$ in $C([\varepsilon, 1],
\mathcal{S}'(\R^d))$ as $n\to\infty$, where $\zeta$ is a standard
normal random variable and
$$C_1=\sqrt{\frac{2\gamma\kappa}{\theta(2\pi)^d}\int_{\R^d}\frac{1}
 {(1+\sum_{k=1}^d|y_k|^{\alpha_k})^3}\d
 y}.$$
\end{thm}

\begin{thm}\label{s2-thm-2} When $1<\bar{\alpha}<2$, let
$F_n^2=n^{(3-\bar{\alpha})\kappa}$. Then for any $\varepsilon>0$, $X_n\Rightarrow
C_2\lambda\zeta$ in $C([\varepsilon, 1], \mathcal{S}'(\R^d))$ as
$n\to\infty$, where $\zeta$ is a standard normal random variable and
$$C_2=\sqrt{\frac{\gamma}{\pi^d}\prod_{k=1}^d\frac{\Gamma(1/\alpha_k)}{\alpha_k}\int_0^\infty\e^{-\theta u}\d u
 \int_0^\infty\e^{-\theta v}\d v\int_0^{v\wedge u}\frac{\e^{\theta s}\d s}{(u+v-2s)^{\bar{\alpha}}}}.$$
\end{thm}

\begin{rem}
(1)  Compared with the corresponding results in the case of weak
degeneration (see \cite[Theorem 2.1, Theorem 2.2 and Remark
2.1]{LX10}), the limit processes are simpler in the temporal structure.

(2)Though Li and Xiao \cite{LX10} pointed out that their results under the case of $\bar{\alpha}=2$ can be strengthened to the weak functional convergence in $C([0, 1],\mathcal{S}'(\R^d))$ by a lengthy and tedious method, due to the strong degeneration, in this paper, we can use a relatively simple way to get the weak convergence in $C([\varepsilon, 1],\mathcal{S}'(\R^d))$.

(3) Note that $X_n(\cdot)\in C([0, 1],
\mathcal{S}'(\R^d))$ and $X_n(0)=0$ and that $X$ is a non-zero time-independent measure-valued Wienner processes.  $X_n$ does not weakly converge to $X$ in $C([0, 1],
\mathcal{S}'(\R^d))$ because  the limit of $X_n$ in $C([0, 1],
\mathcal{S}'(\R^d))$ must be continuous and its initial value has to be $0$ a.s.

(4) For the case of large dimensions a similar result holds, i.e., if $\bar{\alpha}>2$ and $F_n^2=n^{\kappa}$,
then, for any $\varepsilon>0$, $X_n\Rightarrow X$ in $C([\varepsilon, 1],
\mathcal{S}'(\R^d))$ as $n\to\infty$, where $X$ is a centered
time-independent Gaussian process valued in $\mathcal{S}'(\R^d)$, with
covariance function
 \beqnn
 &&{\rm Cov}(\langle X(s),\phi_1\rangle,\langle X(t), \phi_2\rangle)
 \\&&\qquad=\frac{1}{\theta(2\pi)^d}\int_{\R^d}
\Big[\frac{2}{\sum_{k=1}^d|z_k|^{\alpha_k}}+\frac{\gamma}{(\sum_{k=1}^d|z_k|^{\alpha_k})^2}\Big]
\widehat{\phi}_1(z)\overline{\widehat{\phi}_2(z)}\d
 z.
 \eeqnn
It is also interesting to study the properties of the limit processes. We will discuss these problems elsewhere.
\end{rem}

 For the convenience of
reference, at the end of this section, we collect some formulas and results as
follows.
\begin{lem}(\cite[Remark 2.3]{LX10})\label{s2-lm-1}
Let $z=(z_1,\cdots, z_d)$.  For any $\{\alpha_k>0, k=1, \cdots,
 d\}$, if $0<r<\bar{\alpha}$, then
$ \int_{[0, 1]^d}\frac{1}{\sum_{k=1}^d|z_k|^{r\alpha_k}}\d
z<\infty, $ and if $r>\bar{\alpha}$, then
 $ \int_{\R^d\setminus[0,
1]^d}\frac{1}{\sum_{k=1}^d|z_k|^{r\alpha_k}}\d z<\infty.
 $
 Therefore, if $\tau(z)$
is bounded and $\int_{\R^d} \tau(z)\d z<\infty$, then
  $ \int_{\R^d}\frac{\tau(z)}{\sum_{k=1}^d|z_k|^{r\alpha_k}}\d
 z<\infty,$
for all $r\in(0, \bar{\alpha})$.
\end{lem}

Let $\phi_1$, $\phi_2$ and $\phi_3$
be functions from $\R^d$ to $\R$, bounded and integrable. Then
 \beqlb\label{s2-3}
 \int_{\R^d}\phi_1(x)\phi_2(x)\d x&=& \frac{1}{(2\pi)^{d}}\int_{\R^d}\widehat{\phi}_1(z)\overline{\widehat{\phi}_2(z)}\d
 z,
 \eeqlb
(the Plancherel formula). Furthermore, if $\widehat{\phi_1}$ and $\widehat{\phi_2}$  are integrable, then
 \beqlb\label{s2-4}
 \int_{\R^d}\phi_1(x)\phi_2(x)\phi_3(x)\d x&=&\frac{1}{(2\pi)^{2d}}\int_{\R^{2d}}
 \widehat{\phi}_1(z)\widehat{\phi}_2(z_1)\overline{\widehat{\phi}_3(z+z_1)}\d
 z\d z_1,
 \eeqlb
(the inverse Fourier transform), and that, by the
 Riemann-Lebesgue Lemma, $\widehat{\phi}_1(z)$ is bounded and goes to $0$ as $|z|\to\infty$.

Since components of
$\vec{\xi}$ are symmetric stable L\'{e}vy processes and independent
of each other, for any $t>0$
 \beqlb\label{s2-5}
 \int_{\R^d} \phi_1(x)T_t\phi_2(x)\d x&=&\int_{\R^d}\phi_2(x)T_t\phi_1(x)\d
x,
 \eeqlb
and
 \beqlb\label{s2-6}
 \widehat{T_t\phi}_1(z)=\widehat{\phi}_1(z)\e^{-t\sum_{k=1}^d|z_k|^{\alpha_k}}.
 \eeqlb

\section{The proofs of main results}
First of all, we define a sequence of random variables $\tilde{X}_n$ in
$\mathcal{S}'(\R^{d+1})$ as follows:

For any $n\geq0$ and
$\psi\in\mathcal{S}(\R^{d+1})$, let
 \beqlb\label{s3-1}
 \langle\tilde{X}_n, \psi\rangle=\int_0^1\langle X_n(t), \psi(\cdot, t)\rangle\d t.
 \eeqlb

In order to prove the main results, as what Bojdecki et al did in
their serial papers (\cite{BGT061}-\cite{BGT072}), we need show the following facts.

 (i) $\langle\tilde{X}_n, \psi\rangle$ converges in distribution to
 $\langle \tilde{X},\psi\rangle$ for all $\psi\in\mathcal{S}(\R^{d+1})$ as
 $n\to\infty$, where $\tilde{X}_n$ and $\tilde{X}$ are defined as
 (\ref{s3-1}) and $X$ is the corresponding limit process.

 (ii) For any given $\varepsilon>0$, $\{\langle X_n, \phi\rangle, n\geq 1\}$ is tight in $C([\varepsilon, 1])$ for all $\phi\in\mathcal{S}(\R^{d})$, where the
theorem of Mitoma \cite{M83} is used.

 As explained in Bojidecki et al \cite{BGT061}, (i) will be proved if we show that
 \beqlb\label{s4-1}
 \lim_{n\to\infty}\bfE(\e^{-\langle\tilde{X}_n, \psi
 \rangle})=\exp\Big\{\frac{1}{2}\int_0^1\int_0^1 Cov(\langle X(s), \psi(\cdot, s)\rangle,\langle X(t),\psi(\cdot, t)\rangle)\d s\d
 t\Big\},\qquad
 \eeqlb
for each non-negative $\psi\in\mathcal{S}(\R^{d+1})$.

Below, we state the proof of Theorem \ref{s2-thm-1} in detail. Since the proof of Theorem \ref{s2-thm-2} is similar and easier, we omit it.

\medskip

\noindent{\bf Proof of Theorem \ref{s2-thm-1}}. To prove
(\ref{s4-1}), we assume $\psi(x, t)=\phi(x)h(t)$, where
$\phi\in\mathcal{S}(\R^{d})$ and $h\in\mathcal{S}(\R)$ are arbitrary given
nonnegative functions. For general $\psi$, the proof is the same
with slightly more complicated notation.

Now, we recall some formulas from Li and Xiao \cite{LX10} as follows.
 \beqlb\label{s3-2}
   \bfE(\e^{-\langle \tilde{X}_n, \psi\rangle})&=&\exp\Big\{\int_{\R^d}\d x\int_0^n f_n(s)T_s\psi_n(x, s)\d s
   -\int_{\R^d}[1-H_{n,\psi_n}(x,n,0)]\d x\Big\}\nonumber
 \\&=&\exp\Big\{\int_{\R^d} [J_{n,\psi_n}(x, n,
 0)- V_{n,\psi_n}(x, n, 0)]\d x\Big\}\nonumber
 \\&=&\exp\big\{I_1(n,\psi_n)+I_2(n,\psi_n)+I_3(n,\psi_n)\big\},
   \eeqlb
 where
   \beqlb
  I_1(n,\psi_n)&=&\frac{\gamma}{2}\int_{\R^d}\d x\int_0^n\e^{-\delta_n s}V_{n,\psi_n}^2(x, n-s, s)\d
 s,\label{s3-3}
 \\ I_2(n,\psi_n)&=&\int_{\R^d}\d x\int_0^n\e^{-\delta_n s}\psi_n(x,
  s)V_{n,\psi_n}(x, n-s, s)\d s,\label{s3-4}
 \\ I_3(n,\psi_n)&=&\delta_n\int_{\R^d}\d x\int_0^n\e^{-\delta_n s}\d s\int_0^{n-s}\e^{-\gamma u}\chi_{n,\psi_n}(x, u, s) \d
  u.\label{s3-5}
  \eeqlb
 Here
 \beqlb
 \psi_n(x,
s)&=&\frac{1}{F_n}\phi(x)\tilde{h}(\frac{s}{n})\;\;\;\text{and}\;\;\;
\tilde{h}(s)=\int_s^1 h(t)\d t,\label{s3-6}
\\ V_{n,\psi_n}(x, t,
r)&=&1-\bfE_x\Big(\exp\Big\{-\int_0^t\big\langle N_n(s),
\psi_n(\cdot, r+s)\big\rangle\d
 s\Big\}\Big),\label{s3-7}
 \eeqlb
and
 \beqlb\label{s3-8}
 \chi_{n,\psi_n}(x, u,
s)=\bfE_x\Big[\Big(1-\e^{-\int_0^u\psi_n(\vec{\xi}_n(v), s+v)\d
v}\Big)\psi_n(\vec{\xi}_n(u), s+u)\Big].
 \eeqlb
In addition, for any $x\in\R^d$ and $t, s\geq 0$, from Li and Xiao
\cite{LX10}  we still have that
 \beqlb\label{s3-9}
  V_{n,\psi_n}(x, t,
 r)\leq\int_0^t f_n(s)T_s\psi_n(\cdot, r+s)(x)\d s=:J_{n,\psi_n}(x, t,
 r),
 \eeqlb
and that
 \beqlb\label{s3-10}
 J_{n,\psi_n}(x, t,
 r)- V_{n,\psi_n}(x, t,
 r)&=&\delta_n\int_0^t\e^{-\delta_n s}T_s\Big(\int_0^{t-s}\e^{-\gamma u}
 \chi_n(\cdot, u, r+s)
 \d u\Big)(x)\d s\nonumber
 \\&&+\int_0^t\e^{-\delta_n s}T_s\Big\{\psi_n(\cdot, r+s)V_{n,\psi_n}(\cdot, t-s,
 r+s)\Big\}(x)\d s\nonumber
 \\&&+\frac{\gamma}{2}\int_0^t\e^{-\delta_n s}T_sV_{n,\psi_n}^2(x, t-s, r+s)\d
 s.\qquad
 \eeqlb

Below, we discuss the limits of $I_1(n,\psi_n)$, $I_2(n,\psi_n)$ and $I_3(n,\psi_n)$, respectively. We remind that for all $t>0$ and $y=t^Hz$, where $H$ be the $d\times d$ diagonal matrix $(1/\alpha_k)_{1\leq
k\leq d}$, $$\d y=t^2\d z.$$

\medskip

{\bf Step 1} We are going to get the limit of $I_1(n,\psi_n)$. From (\ref{s3-3}), we get that
  \beqlb\label{s3-3-1}
  I_1(n,\psi_n)=I_{11}(n,\psi_n)+I_{12}(n,\psi_n),
  \eeqlb
where
 \beqlb
 I_{11}(n,\psi_n)&=&\frac{\gamma}{2}\int_{\R^d}\d x\int_0^n\e^{-\delta_n
 s}J_{n,\psi_n}^2(x, n-s, s)\d s,\label{s3-3-2}
 \\I_{12}(n,\psi_n)&=&\frac{\gamma}{2}\int_{\R^d}\d x\int_0^n\e^{-\delta_n
 s}(V_{n,\psi_n}^2(x, n-s, s)-J_{n,\psi_n}^2(x, n-s, s))\d
 s.\qquad\label{s3-3-3}
 \eeqlb

We first consider the limit of $I_{11}(n,\psi_n)$. Substituting (\ref{s1-2}), (\ref{s3-6}) and (\ref{s3-9}) into (\ref{s3-3-2}),  we get that
  \beqnn
  I_{11}(n,\psi_n)&=&\frac{n^3\gamma2}{2F_n^2}\int_0^1\e^{-n\delta_n s}\d s\int_{\R^d}\Big[\int_0^{1-s}\bar{f}_n(nu)\e^{-n\delta_n u}\tilde{h}(s+u)T_{nu}\phi(x)\d
 u\Big]^2\d x.
 \eeqnn
Furthermore applying (\ref{s2-3}),(\ref{s2-5}) and (\ref{s2-6}) to the above formula, and noting that
$\bar{f}_n(u)$ converges uniformly to $1$ as $n\to\infty$, we derive that
\beqlb\label{s3-31-1}
 \lim_{n\to\infty}I_{11}(n,\psi_n)&=&\lim_{n\to\infty}\frac{n^3\gamma}{2F_n^2}
 \int_0^1\e^{-n\delta_n s}\d s\int_{\R^d}\Big[\int_0^{1-s}\e^{-n\delta_n u}\tilde{h}(s+u)T_{nu}\phi(x)\d
 u\Big]^2\d x\nonumber
\\&=&\lim_{n\to\infty}\bigg\{\frac{ n^3\gamma/2}{(2\pi)^dF_n^2}
 \int_0^1\e^{-n\delta_n s}\d s\int_{\R^d}|\widehat{\phi}(z)|^2\nonumber
 \\&&\qquad\quad\times\Big[\int_s^1\e^{-n(u-s)(\delta_n+\sum_{k=1}^d|z_k|^{\alpha_k})}\tilde{h}(u)\d
 u\Big]^2\d z\bigg\}.
 \eeqlb
Substituting $\tilde{h}(u)=\int_u^1h(t)\d t$ and
$F_n^2=n^{\kappa}\ln n$ into (\ref{s3-31-1}), by changing the
integral order we obtain that
 \beqlb\label{s3-31-2}
 &&\lim_{n\to\infty}I_{11}(n,\psi_n)=\lim_{n\to\infty}\Big\{\frac{n\gamma/2}{(2\pi)^dF_n^2}
 \int_0^1 h(r)\d r\int_0^1 h(t) \d t\int_0^{t\wedge r}\e^{-n\delta_n s}\d
 s\int_{\R^d}|\widehat{\phi}(z)|^2\nonumber
 \\&&\qquad\qquad\qquad\qquad\times\frac{(1-\e^{-n(r-s)(\delta_n+\sum_{k=1}^d|z_k|^{\alpha_k})})
 (1-\e^{-n(t-s)(\delta_n+\sum_{k=1}^d|z_k|^{\alpha_k})})}
 {(\delta_n+\sum_{k=1}^d|z_k|^{\alpha_k})^2}\d z\Big\}\nonumber
 \\&&\qquad\qquad=\lim_{n\to\infty}\frac{\gamma}{(2\pi)^d}\int_0^1 h(r)\d r\int_0^r
 h(t)\d t\int_0^{\infty}\e^{-n^\kappa\delta_n s}\frac{W_{n, r, t, s}(n^{\kappa}\delta_n)}{\ln n}\d
 s,\qquad\quad
 \eeqlb
where for any $x>0$, $0\leq t\leq r\leq 1$ and $s>0$,
 \beqnn
 W_{n, r, t, s}(x)=\begin{cases}\int_{\R^d}|\widehat{\phi}(z)|^2\Phi(x, s, n^{1-\kappa}
r, n^{1-\kappa} t, n^{\kappa}\sum_{k=1}^d|z_k|^{\alpha_k})n^{2\kappa}\d z, & s<n^{1-\kappa}t;
\\ 0,&s\geq n^{1-\kappa}t,
\end{cases}
\eeqnn
and for any $x>0$, $s, u, v, y\in[0,\infty)$
\beqnn
\Phi(x, s, u,v,y)=\frac{(1-\e^{-(u-s)(x+y)})
 (1-\e^{-(v-s)(x+y)})}
 {(x+y)^2}.
 \eeqnn
Let
 \beqlb\label{s3-31-6}
 \tilde{W}(n)=\int_{\R^d}|\widehat{\phi}(z)|^2\frac{(1-\e^{-
 n\sum_{k=1}^d|z_k|^{\alpha_k}})^2}
 {(\sum_{k=1}^d|z_k|^{\alpha_k})^2\ln n}\d z.
 \eeqlb
It is easy to see that
$$W_{n, r, t, s}(n^\kappa\delta_n)/\ln n\leq \tilde{W}(n)$$
 for all $(r, t,
s)\in\{0\leq t\leq r\leq 1; 0\leq s\}$, where we use the decreasing of $(1-\e^{-r})/r$ on $r\in(0,+\infty)$. Furthermore, applying L'H\^{o}pital's law, we obtain that
 \beqlb\label{s3-31-61}
 \lim_{n\to\infty}\tilde{W}(n)&=&\lim_{n\to\infty}\int_{\R^d}2n|\widehat{\phi}(z)|^2\frac{(1-\e^{-
 n\sum_{k=1}^d|z_k|^{\alpha_k}})\e^{-
 n\sum_{k=1}^d|z_k|^{\alpha_k}}} {\sum_{k=1}^d|z_k|^{\alpha_k}}\d z\nonumber
 \\&=&2|\widehat{\phi}(0)|^2\int_{\R^d}\frac{\e^{-\sum_{k=1}^d|z_k|^{\alpha_k}}
 (1-\e^{-\sum_{k=1}^d|z_k|^{\alpha_k}})}{\sum_{k=1}^d|z_k|^{\alpha_k}}\d z<\infty,
 \eeqlb
and hence $\{\tilde{W}(n)\}$ is bounded. Therefore, the dominated convergence theorem plus the convergence of $n^\kappa\delta_n\to\theta$
yields that if
\beqlb\label{s3-31-8}
W_{n, r, t, s}(n^\kappa\delta_n)/\ln n\to 2\kappa\int_{\R^d}\frac{|\widehat{\phi}(0)|^2}{(1+\sum_{k=1}^d|y_k|^{\alpha_k})^3}\d
 y,\;\; a.s.
\eeqlb
on $(r, t,
s)\in\{0\leq t\leq r\leq 1, 0\leq s\}$, then
 \beqlb\label{s3-31-7}
 I_{11}(n,\psi_n)&\to&\frac{2\gamma\kappa}{(2\pi)^d}
 \int_{\R^d}\frac{|\widehat{\phi}(0)|^2}
 {(1+\sum_{k=1}^d|y_k|^{\alpha_k})^3}\d
 y\int_0^1 h(r)\d r\int_0^rh(t)\d t\int_0^\infty\e^{-\theta s}\d s\nonumber
 \\&=&\frac{\gamma\kappa}{\theta(2\pi)^d}\int_{\R^d}\frac{\d
 y}
 {(1+\sum_{k=1}^d|y_k|^{\alpha_k})^3}\Big(\int_{\R^d} \phi(x)\d x\int_0^1 h(t)\d t\Big)^2.\qquad
 \eeqlb
Below, we prove (\ref{s3-31-8}). To this end,  by the mean-value theorem and using the substitution $y=\Theta_n z:=(n^{\kappa})^H z$ we have that
\beqlb
 &&|W_{n, r, t, s}(x_1)-W_{n, r, t,
s}(x_2)|\nonumber
\\&&=\begin{cases}|x_1-x_2|\int_{\R^d}|\widehat{\phi}(\Theta_n^{-1}y)|^2\Phi'_x(\vartheta, s, n^{1-\kappa}
r, n^{1-\kappa} t, \sum_{k=1}^d|y_k|^{\alpha_k})\d y, & s<n^{1-\kappa}t,
\\ 0,&s\geq n^{1-\kappa}t,
\end{cases}\nonumber
\eeqlb
for any $x_1, x_2>0$, where $\vartheta\in(x_1, x_2)$. Note that
 \beqnn
 &&\Phi'_x(\vartheta, s, n^{1-\kappa}
r, n^{1-\kappa} t, \sum_{k=1}^d|y_k|^{\alpha_k})
\\&&=\frac{(n^{1-\kappa} r-s)\e^{-(n^{1-\kappa }r-s)(\vartheta+\sum_{k=1}^d|y_k|^{\alpha_k})}(1-\e^{-(n^{1-\kappa }t-s)(\vartheta+\sum_{k=1}^d|y_k|^{\alpha_k})})}{(\vartheta+\sum_{k=1}^d|y_k|^{\alpha_k})^2}
 \\&&\quad+\frac{(n^{1-\kappa} t-s)\e^{-(n^{1-\kappa }t-s)(\vartheta+\sum_{k=1}^d|y_k|^{\alpha_k})}(1-\e^{-(n^{1-\kappa }r-s)(\vartheta+\sum_{k=1}^d|y_k|^{\alpha_k})})}{(\vartheta+\sum_{k=1}^d|y_k|^{\alpha_k})^2}
 \\&&\quad-\frac{2(1-\e^{-(n^{1-\kappa }r-s)(\vartheta+\sum_{k=1}^d|y_k|^{\alpha_k})})(1-\e^{-(n^{1-\kappa }t-s)(\vartheta+\sum_{k=1}^d|y_k|^{\alpha_k})})}{(\vartheta+\sum_{k=1}^d|y_k|^{\alpha_k})^3}.
 \eeqnn
 We have that for any given $1\geq r\geq t\geq 0$ and $s\geq0$, there exists $N>0$ such that for all $n>N$,
 \beqlb\label{s3-31-3}
 |W_{n, r, t, s}(x_1)-W_{n, r, t,
s}(x_2)|\leq|x_1-x_2|\int_{\R^d}|\widehat{\phi}(\Theta_n^{-1}y)|^2Z_n(\vartheta,
r, t, s, y) \d y,\qquad
 \eeqlb
 where
 \beqnn
 Z_n(\vartheta, r, t, s, y)&=&\frac{2}{(\vartheta+\sum_{k=1}^d|y_k|^{\alpha_k})^3}+\frac{(n^{1-\kappa} r-s)\e^{-(n^{1-\kappa} r-s)\sum_{k=1}^d|y_k|^{\alpha_k}}}{\vartheta^2}
 \\&&+\frac{(n^{1-\kappa} t-s)\e^{-(n^{1-\kappa} t-s)(\sum_{k=1}^d|y_k|^{\alpha_k})}}{\vartheta^2}.
 \eeqnn
Because for any given $r, t, s$ and sufficiently large $n$, $\int_{\R^d}Z_n(\vartheta,r,
t, s, y)\d y$ equals
 \beqnn
 \frac{1}{\vartheta^2}\Big[\frac{1}{n^{1-\kappa} r-s}+\frac{1}{n^{1-\kappa} t-s}\Big]
 \int_{\R^d}\e^{-\sum_{k=1}^d|y_k|^{\alpha_k}}\d y+\int_{\R^d}\frac{1}
 {(\vartheta+\sum_{k=1}^d|y_k|^{\alpha_k})^3}\d y,
 \eeqnn
which is bounded for sufficiently large $n$, from
$n^{1-\kappa}\delta_n\to\theta\in(0,+\infty)$ and (\ref{s3-31-3}),
we obtain that as $n\to\infty$,
 \beqlb\label{s3-31-4}
 |W_{n, r, t, s}(n^{1-\kappa}\delta_n)-W_{n, r, t, s}(\theta)|\to
 0,
 \eeqlb
for any given $(r, t, s)\in\{ 0\leq t\leq r\leq 1, 0\leq
s\}$. Therefore
 \beqlb\label{s3-31-5}
 \lim_{n\to\infty}\frac{W_{n, r, t, s}(n^{1-\kappa}\delta_n)}{\ln n}=\lim_{n\to\infty}\frac{W_{n, r, t,
 s}(\theta)}{\ln n}=\lim_{n\to\infty} n\frac{\partial W_{n, r, t,
 s}(\theta) }{\partial n},\qquad
 \eeqlb
where we use L'H\^{o}pital's law at the second equality. Note that for any $(r, t, s)\in\{ 0\leq t\leq r\leq 1, 0\leq
s\}$ and sufficiently large $n$
 \beqnn
n\frac{\partial W_{n, r, t, s}(\theta)
}{\partial
n}&=&n\int_{\R^d}|\widehat{\phi}(z)|^2\frac{\partial\Phi(\theta, s, n^{1-\kappa}
r, n^{1-\kappa} t, n^{\kappa}\sum_{k=1}^d|z_k|^{\alpha_k})n^{2\kappa}}{\partial n}\d z,
\eeqnn
which, by direct calculations, equals
 \beqnn
&&\int_{\R^d}\Big\{|\widehat{\phi}(z)|^2\big((1-\kappa)\theta
n^{1-\kappa} r+(n^{1-\kappa} r-\kappa s)n^{\kappa}\sum_{k=1}^d|z_k|^{\alpha_k}\big)
\\&&\times\frac{\e^{-(n^{1-\kappa}r-s)(\theta+n^\kappa\sum_{k=1}^d|z_k|^{\alpha_k})}
(1-\e^{-(n^{1-\kappa}t-s)(\theta+n^\kappa\sum_{k=1}^d|z_k|^{\alpha_k})})}
{(\theta+n^\kappa\sum_{k=1}^d|z_k|^{\alpha_k})^2} \Big\}n^{2\kappa}\d
z\\&+&\int_{\R^d}\Big\{|\widehat{\phi}(z)|^2\big((1-\kappa)\theta
n^{1-\kappa} t+(n^{1-\kappa} t-\kappa s)n^{\kappa}\sum_{k=1}^d|z_k|^{\alpha_k}\big)
\\&&\times\frac{\e^{-(n^{1-\kappa}t-s)(\theta+n^\kappa\sum_{k=1}^d|z_k|^{\alpha_k})}
(1-\e^{-(n^{1-\kappa}r-s)(\theta+n^\kappa\sum_{k=1}^d|z_k|^{\alpha_k})})}
{(\theta+n^\kappa\sum_{k=1}^d|z_k|^{\alpha_k})^2} \Big\}n^{2\kappa}\d
z
 \\&+&\int_{\R^d}\Big\{|\widehat{\phi}(z)|^2\frac{2\theta\kappa(1-\e^{-(n^{1-\kappa}r-s)(\theta+n^\kappa\sum_{k=1}^d|z_k|^{\alpha_k})})
 } {(\theta+n^\kappa\sum_{k=1}^d|z_k|^{\alpha_k})^3}
 \\&&\times(1-\e^{-(n^{1-\kappa}t-s)(\theta+n^\kappa\sum_{k=1}^d|z_k|^{\alpha_k})})\Big\}n^{2\kappa}\d
 z.
 \eeqnn
Now substituting $y=(n^\kappa)^H z$ into the above formula and letting $n\to\infty$, we get that
for any given $r, t, s$, as $n\to\infty$,
 \beqnn
 n\frac{\partial W_{n, r, t,
 s}(\theta) }{\partial
 n}&\to&
 2\kappa|\widehat{\phi}(0)|^2\int_{\R^d}\frac{\theta}{(\theta+\sum_{k=1}^d|y_k|^{\alpha_k})^3}\d
 y\nonumber
 \\&=&2\kappa|\widehat{\phi}(0)|^2\int_{\R^d}\frac{1}{(1+\sum_{k=1}^d|y_k|^{\alpha_k})^3}\d
 y,
 \eeqnn
which is the desired formula (\ref{s3-31-8}).

 To study the limit of $I_{12}(n,\psi_n)$, we first observe that  from
(\ref{s3-8})-(\ref{s3-10}),
 \beqlb\label{s3-3-13}
  &&J^2_{n,\psi_n}(x, n-s, s)-V^2_{n,\psi_n}(x, n-s, s)\nonumber
 \\&&\qquad\leq 2\bigg[\delta_n\int_0^{n-s}\e^{-\delta_n u}T_u\Big(\int_0^{n-s-u}\e^{-\gamma v}\chi_{n,\psi_n}(\cdot, v, s+u)\d v\Big)(x)\d u\nonumber
 \\&&\qquad+\int_0^{n-s}\e^{-\delta_n u}T_u\big(\psi_n(\cdot, s+u)J_{n,\psi_n}(\cdot, n-s-u,
 s+u)\big)(x)\d u\nonumber
 \\&&\qquad+\frac{\gamma}{2}\int_0^{n-s}\e^{-\delta_n u}T_uJ_{n,\psi_n}^2(x, n-s-u, s+u)\d
 u\bigg]J_{n,\psi_n}(x, n-s, s)\nonumber
 \\&&\qquad=:2J_{n,\psi_n}(x, n-s, s)\Big(\bar{I}_{121}(n,\psi_n)+\bar{I}_{122}(n,\psi_n)+\bar{I}_{123}(n,\psi_n)\Big).
 \eeqlb
Let
 \beqlb
 I_{121}(n,\psi_n)&:=&\gamma\int_{\R^d}\d x\int_0^n\e^{-\delta_n
 s}J_{n,\psi_n}(x, n-s, s)\bar{I}_{121}(n,\psi_n)\d s,\label{s3-3-14}
\\I_{122}(n,\psi_n)&:=&\gamma\int_{\R^d}\d x\int_0^n\e^{-\delta_n
 s}J_{n,\psi_n}(x, n-s, s)\bar{I}_{122}(n,\psi_n)\d s,\label{s3-3-15}
 \\I_{123}(n,\psi_n)&:=&\gamma\int_{\R^d}\d x\int_0^n\e^{-\delta_n
 s}J_{n,\psi_n}(x, n-s, s)\bar{I}_{123}(n,\psi_n)\d s. \label{s3-3-16}
 \eeqlb
From Li and Xiao \cite[(3.45), (3.47) and (3.48)]{LX10}  we derive that
 \beqlb
 I_{121}(n,\psi_n)&\leq&
 \frac{2\delta_n}{\gamma(\gamma-\delta_n)}I_{11}(n,\psi_n)\label{s3-3-17},
 \eeqlb
and that
 \beqlb
 I_{122}(n,\psi_n)&\leq&\frac{Kn}{F_n^3} \int_0^1\e^{-n\delta_n s}\d s\int_{\R^{2d}}\frac{|\widehat{\phi}(z)|(1-\e^{-n\sum_{k=1}^d|z_k|^{\alpha_k}})^2}{\big(\sum_{k=1}^d|z_k|^{\alpha_k}\big)^2}
 \frac{|\widehat{\phi}(z')|\d z\d
 z'}{\sum_{k=1}^d|z'_k|^{\alpha_k}};\qquad\label{s3-3-18}
\\  I_{123}(n,\psi_n)&\leq&\frac{Kn}{F_n^3}\int_0^1\e^{-n\delta_n s}\d s\int_{\R^{2d}}\frac{1-\e^{-n\sum_{k=1}^d|z_k+z'_k|^{\alpha_k}}}
 {\sum_{k=1}^d|z_k+z'_k|^{\alpha_k}}
 \frac{1-\e^{-n\sum_{k=1}^d|z_k|^{\alpha_k}}}{\sum_{k=1}^d|z'_k|^{\alpha_k}}
 \nonumber
 \\&&\qquad\qquad\times\frac{(1-\e^{-n\sum_{k=1}^d|z_k|^{\alpha_k}})^2}{(\sum_{k=1}^d|z_k|^{\alpha_k})^2}
 |\widehat{\phi}(z)||\widehat{\phi}(z')||\widehat{\phi}(z+z')|\d z\d
 z'.\label{s3-3-19}
 \eeqlb

Firstly, since $\delta_n\to 0$, from (\ref{s3-31-7}) and (\ref{s3-3-17}) it follows that
 \beqlb\label{s3-3-20}
 I_{121}(n,\psi_n)\to0.
 \eeqlb

Secondly, applying the fact $F_n^2=n^{\kappa}\ln n$ to
(\ref{s3-3-18})  we get that
   \beqlb\label{s3-3-21}
 I_{122}(n,\psi_n)\leq\frac{K}{n^{\kappa}\delta_n F_n}\int_{\R^{2d}}
 \frac{|\widehat{\phi}(z)|(1-\e^{-n\sum_{k=1}^d|z_k|^{\alpha_k}})^2}{\ln n\big(\sum_{k=1}^d|z_k|^{\alpha_k}\big)^2}
 \frac{|\widehat{\phi}(z')|}{\sum_{k=1}^d|z'_k|^{\alpha_k}}\d z\d
 z'.\qquad
 \eeqlb
Since $n^{\kappa}\delta_n\to\theta\in(0,\infty)$, from (\ref{s3-31-6}), (\ref{s3-31-61}),
(\ref{s3-3-21})and Lemma \ref{s2-lm-1} we derive that
 \beqlb\label{s3-3-22}
 I_{122}(n,\psi_n)\to 0.
 \eeqlb

At last, using the fact $F_n^2=n^{\kappa}\ln n$ again, from
 (\ref{s3-3-19})
 we  get that for some constant $K>0$,
 \beqnn
  I_{123}(n,\psi_n)&\leq&\frac{K}{n^{\kappa}\delta_n
 n^{\frac{\kappa}{2}}(\ln n)^{3/2}}\int_{\R^{2d}}\frac{1-\e^{-n\sum_{k=1}^d|z_k+z'_k|^{\alpha_k}}}
 {\sum_{k=1}^d|z_k+z'_k|^{\alpha_k}}
 \frac{1-\e^{-n\sum_{k=1}^d|z_k|^{\alpha_k}}}{\sum_{k=1}^d|z'_k|^{\alpha_k}}
 \nonumber
 \\&&\qquad\qquad\qquad\times\frac{(1-\e^{-n\sum_{k=1}^d|z_k|^{\alpha_k}})^2}{(\sum_{k=1}^d|z_k|^{\alpha_k})^2}
 |\widehat{\phi}(z)||\widehat{\phi}(z')||\widehat{\phi}(z+z')|\d z\d
 z'.\qquad
 \eeqnn
 Furthermore, by using the  inequality $1-\e^{-x}\leq
 x^{\kappa/8}$ for $x\geq 0$ we have that
  \beqlb\label{s3-3-23}
  I_{123}(n,\psi_n)&\leq&\frac{K}{n^{\kappa}\delta_n
 (\ln n)^{3/2}}\int_{\R^{2d}}\frac{|\widehat{\phi}(z+z')|}
 {(\sum_{k=1}^d|z_k+z'_k|^{\alpha_k})^{1-\kappa/8}}
 \frac{|\widehat{\phi}(z')|}{(\sum_{k=1}^d|z'_k|^{\alpha_k})^{1-\kappa/8}}
 \nonumber
 \\&&\qquad\qquad\qquad\times\frac{|\widehat{\phi}(z)|}{(\sum_{k=1}^d|z_k|^{\alpha_k})^{2-\kappa/4}}
 \d z\d
 z'.
  \eeqlb
 Since $\bar{\alpha}=2$, from Lemma \ref{s2-lm-1} we know
 $$\int_{\R^d}\frac{|\widehat{\phi}(z)|^2}{(\sum_{k=1}^d|z_k|^{\alpha_k})^{2-\kappa/4}}\d z<\infty,$$
and hence by H\"{o}lder inequality,
  $$\int_{\R^d}\frac{|\widehat{\phi}(z')||\widehat{\phi}(z+z')|} {(\sum_{k=1}^d|z_k+z'_k|^{\alpha_k})^{1-\kappa/8}(\sum_{k=1}^d|z'_k|^{\alpha_k})^{1-\kappa/8}}\d
 z'$$
is bounded for all $z\in\R^d$. Therefore
 $$\int_{\R^{2d}}\frac{|\widehat{\phi}(z+z')|}{(\sum_{k=1}^d|z_k+z'_k|^{\alpha_k})^{1-\kappa/8}}
 \frac{|\widehat{\phi}(z')|}{(\sum_{k=1}^d|z'_k|^{\alpha_k})^{1-\kappa/8}}
 \frac{|\widehat{\phi}(z)|}{(\sum_{k=1}^d|z_k|^{\alpha_k})^{2-\kappa/4}}\d z\d
 z'<\infty.$$
Hence, (\ref{s3-3-23}) and the fact that
$n^{\kappa}\delta_n\to\theta\in(0,\infty)$ imply that as $n\to\infty$,
 \beqlb\label{s3-3-24}
 I_{123}(n,\psi_n)\to 0.
 \eeqlb
Consequently, from (\ref{s3-3-13})-(\ref{s3-3-20}), (\ref{s3-3-22}), (\ref{s3-3-24}) and (\ref{s3-3-3})  we have that
 \beqlb\label{s3-3-36}
 I_{12}(n,\psi_n)\to 0.
 \eeqlb
 Combining (\ref{s3-31-7}) and (\ref{s3-3-36}) with (\ref{s3-3-1}) we derive that as $n\to\infty$
 \beqlb\label{s3-3-6}
  I_{1}(n,\psi_n)\to \frac{\gamma\kappa}{\theta(2\pi)^d}\int_{\R^d}\frac{\d
 y}
 {(1+\sum_{k=1}^d|y_k|^{\alpha_k})^3}\Big(\int_{\R^d} \phi(x)\d x\int_0^1 h(t)\d t\Big)^2.
 \eeqlb

{\bf Step 2} We are going to get limits of $I_2(n,\psi_n)$ and $I_3(n,\psi_n)$. Let
 \beqlb\label{s3-4-1}
 \tilde{I}_2(n,\psi_n)&=&\int_{\R^d}\d x\int_0^n\e^{-\delta_n s}\psi_n(x,
 s)J_{n,\psi_n}(x, n-s, s)\d s\nonumber
 \\&=&\int_0^n\e^{-\delta_n s}\d s\int_0^{n-s}f_n(v)\d v\int_{\R^d}\psi_n(x,
 s)T_v\psi_n(x, s+v)\d x.
 \eeqlb
From (\ref{s3-4}) and (\ref{s3-9}), it follows that
 \beqlb\label{s3-4-0}
 I_2(n,\psi_n)\leq \tilde{I}_2(n,\psi_n).\qquad
 \eeqlb
Furthermore, by (\ref{s1-2}),  (\ref{s2-3}),  (\ref{s2-6}) and (\ref{s3-6}), there exists $K>0$ such that
  \beqlb\label{s3-4-2}
 I_2(n,\psi_n)&\leq&\frac{K}{F_n^2}\int_0^n\e^{-\delta_n s}\d s\int_0^{n-s}\d v
 \int_{\R^d}|\widehat{\phi}(z)|^2\e^{-v\sum_{k=1}^d|z_k|^{\alpha_k}}\d
 z\nonumber
 \\&=&\frac{K}{F_n^2}\int_{\R^d}|\widehat{\phi}(z)|^2\d z\int_0^1\e^{-n\delta_n s}n\d s
 \int_0^{1-s}\e^{-nv\sum_{k=1}^d|z_k|^{\alpha_k}}n\d
 v\nonumber
 \\&\leq&\frac{nK}{F_n^2}\int_{\R^d}\frac{|\widehat{\phi}(z)|^2}{\sum_{k=1}^d|z_k|^{\alpha_k}}\d z\int_0^1\e^{-n\delta_n s}\d s.
 \eeqlb
Since $\bar{\alpha}=2$ and $F_n^2=n^{\kappa}\ln n$,
(\ref{s3-4-2}) and Lemma \ref{s2-lm-1} imply that
  \beqlb\label{s3-4-5}
  I_2(n,\psi_n)\leq\frac{K}{\ln n}\int_{\R^d}\frac{|\widehat{\phi}(z)|^2}
  {\sum_{k=1}^d|z_k|^{\alpha_k}}\d
  z\int_0^{n^{1-\kappa}}\e^{-n^{\kappa}\delta_n s}\d s\to 0.
  \eeqlb

To get the limit of $I_3(n,\psi_n)$, we let
 \beqlb\label{s3-2-1}
  \tilde{I}_3(n,\psi_n)&:=&\delta_n\int_0^n\e^{-\delta_n s}\d
 s\int_0^{n-s}\e^{-\gamma u}\d u\nonumber
 \\&&\qquad\quad\times\int_{\R^d}\bfE_x\Big(\int_0^u\psi_n(\vec{\xi}_n(v), s+v)\d v\psi_n(\vec{\xi}_n(u),
 s+u)\Big)\d x.\qquad
 \eeqlb
Then by (\ref{s3-5}), (\ref{s3-8}) and (\ref{s3-2-1}), we have
that
 \beqlb\label{s3-2-0}
 I_3(n,\psi_n)\leq\tilde{I}_3(n,\psi_n),
 \eeqlb
and by (\ref{s2-3}), (\ref{s2-6}), (\ref{s3-6}) and  (\ref{s3-2-1}),  we have that
 \beqlb\label{s3-2-3}
 \tilde{I}_3(n,\psi_n)&\leq&\frac{\delta_n}{F_n^2(2\pi)^d}\int_0^n\e^{-\delta_n s}\d s\int_0^{n-s}\e^{-\gamma u}\d u
  \int_0^u \tilde{h}(\frac{s+v}{n})\tilde{h}(\frac{s+u}{n})\d v\nonumber
  \\&&\qquad\qquad\qquad\qquad\times\int_{\R^d}|\widehat{\phi}(z)|^2
  \e^{-(u-v)\sum_{k=1}^d|z_k|^{\alpha_k}}\d z;
 \eeqlb
see also Li and Xiao \cite[(3.68)]{LX10}. Since $\tilde{h}$ is
bounded and $\int_{\R^d}|\widehat{\phi}(z)|^2\d z<\infty$,
(\ref{s3-2-0}) and (\ref{s3-2-3}) yield that
 $$0\leq I_3(n,\psi_n)\leq K\frac{\delta_n}{F_n^2}\int_0^n\e^{-\delta_n s}\d s
 \int_0^{n-s}\e^{-\gamma u}u\d u\leq\frac{K}{\gamma F_n^2},$$
for some constant $K>0$. Therefore, $F_n\to\infty$ indicates that
 \beqlb\label{s3-2-2}
 I_3(n,\psi_n)\to 0.
 \eeqlb

To get the left hand side of (\ref{s4-1}), we substitute (\ref{s3-3-6}), (\ref{s3-4-5}) and (\ref{s3-2-2}) into (\ref{s3-2}) and obtain that as $n\to\infty$,
 \beqlb
 \bfE(\e^{-\langle \tilde{X}_n, \psi\rangle})&\to&\exp\Big\{\frac{\gamma\kappa}{\theta(2\pi)^d}\int_{\R^d}\frac{\d
 y}
 {(1+\sum_{k=1}^d|y_k|^{\alpha_k})^3}\Big(\int_{\R^d} \phi(x)\d x\int_0^1 h(t)\d t\Big)^2\Big\}\nonumber
 \\&&=\exp\Big\{\frac{1}{2}\int_0^1\int_0^1 {\rm Cov}(\langle X(s), \psi(\cdot, s)\rangle,\langle X(t),\psi(\cdot, t)\rangle)\d s\d
 t\Big\} ,
 \eeqlb
where $X$ is the limit process in Theorem \ref{s2-thm-1}. Therefore (\ref{s4-1}) holds and (i) is proved.


Now we are in the place to prove the tightness of  $\{\langle X_n, \phi\rangle, n\geq 1\}$ in $C([\varepsilon, 1], \R)$. Note that by the same argument as those
used in Bojdecki {\it et al} \cite{BGT072}, we also have that $X_n$ converges to $X$ in finite-dimensional distributions. This implies the tightness of $\{\langle X_n(\varepsilon),\phi\rangle\}$. According to the proof of Proposition 3.3 in \cite{BGT071}, the remainder is
to prove that for all $\phi\in\mathcal{S}(\R^d)$, $\varepsilon\leq
t_1<t_2\leq 1$ and $\eta>0$, there exist constants $a\geq1$, $b>0$
and $K>0$, which is independent of $t_1, t_2$, such that for all
$n\geq1$.
 \beqlb\label{s4-1-1}
 \int_0^{1/\eta}\Big(1-{\rm Re}\Big(\bfE\big(\exp\{-i\omega\langle\tilde{X}_n,
 \phi h\rangle\}\big)\Big)\Big)\d\omega\leq\frac{K}{\eta^{a}}(t_2-t_1)^{1+b},
 \eeqlb
where $h\in\mathcal{S}(R)$ is an approximation of ${\bf
1}_{\{t_2\}}(t)-{\bf 1}_{\{t_1\}}(t)$ supported on $[t_1, t_2]$ such
that
 $\tilde{h}(t)$
satisfies
 \beqlb\label{s3-5-1}
 \tilde{h}\in\mathcal{S}(R),\qquad 0\leq \tilde{h}\leq {\bf 1 }_{[t_1, t_2]}.
 \eeqlb

Repeating the discussion on $\bfE(\exp\{-\langle\tilde{X}_n,
 \psi\rangle\})$ (see Li and Xiao \cite[Section 3]{LX10}) with $\psi$ replaced by $i\omega\phi h$, we can readily  get that
 $$\bfE(\exp\{-i\omega\langle\tilde{X}_n,
 \phi h\rangle\})=\exp\big\{I_1(n, i\omega\psi_n)+I_2(n, i\omega\psi_n)+I_3(n, i\omega\psi_n)\big\},
 $$
and the inequality
 $$
 |V_{n, i\omega\psi_n}|\leq J_{n,\omega\psi_n}=\omega J_{n,\psi_n}.
 $$
Consequently, from the expressions of $I_1, I_2, I_3$ and $I_{11}$,$\tilde{I}_2$, $\tilde{I}_3$ (see (\ref{s3-3})-(\ref{s3-5}), (\ref{s3-3-2}), (\ref{s3-4-1}) and (\ref{s3-2-1}), respectively), it is easy to check that the following inequalities hold.
 \beqlb\label{s4-1-2}
 \begin{cases} |I_1(n, i\omega\psi_n)|\leq I_{11}(n,
 \omega\psi_n)=\omega^2I_{11}(n,\psi_n);
 \\|I_2(n, i\omega\psi_n)|\leq \tilde{I}_2(n,
 \omega\psi_n)=\omega^2\tilde{I}_2(n, \psi_n);
 \\|I_3(n, i\omega\psi_n)|\leq \tilde{I}_3(n,
 \omega\psi_n)=\omega^2\tilde{I}_3(n,\psi_n).
 \end{cases}
 \eeqlb

{\bf (I)} We first estimate the upper bound of $I_{11}(n,\psi_n)$.
Substituting (\ref{s3-6}), (\ref{s3-9}) and (\ref{s1-2}) into
(\ref{s3-3-2}), we get that for some constant $K>0$
 \beqnn
 I_{11}(n,\psi_n)&\leq& \frac{K}{F_n^2}\int_0^n\e^{-\delta_n s}\d
 s\int_0^{n-s}\e^{-\delta_n u}\tilde{h}(\frac{s+u}{n})\d
 u\nonumber
 \\&&\qquad\qquad\qquad\quad\times\int_0^{n-s}\e^{-\delta_n v}\tilde{h}(\frac{s+v}{n})\d v
 \int_{\R^d} T_u\phi(x)T_v\phi(x)\d x.
 \eeqnn
Furthermore, by using (\ref{s2-3}) and (\ref{s2-6}), we have that
\beqnn
 I_{11}(n,\psi_n)&\leq&\frac{2Kn^3}{F_n^2}\int_{\R^d} |\widehat{\phi}(z)|^2\d z\int_0^1\tilde{h}(u)\d
 u\int_0^u\tilde{h}(v)\e^{-n\delta_n (u+v)}\d v\nonumber
 \\&&\qquad\qquad\qquad\quad\times\int_0^v\e^{n\delta_n s}\e^{-n(u+v-2s)\sum_{k=1}^d|z_k|^{\alpha_k}}\d
 s
 \eeqnn
which and the condition (\ref{s3-5-1}) imply that $I_{11}(n,\psi_n)$
is bounded from above by
 \beqnn
 &&\frac{2Kn^3}{F_n^2}\int_{\R^d} |\widehat{\phi}(z)|^2\d
 z\int_{t_1}^{t_2}\d
 u\int_{t_1}^u\e^{-n(u+v)(\delta_n+\sum_{k=1}^d|z_k|^{\alpha_k})}\d v
 \\&&\qquad\qquad\qquad\times\int_0^v\e^{ns(2\sum_{k=1}^d|z_k|^{\alpha_k}+\delta_n)}\d
 s
 \eeqnn
which is further bounded from above by
 \beqlb\label{s3-5-2}
 &&\frac{2Kn^2}{F_n^2}\int_{\R^d} \frac{|\widehat{\phi}(z)|^2}{2\sum_{k=1}^d|z_k|^{\alpha_k}+\delta_n}\d
 z\int_{t_1}^{t_2}\e^{-n\delta_n u}\d
 u\int_{t_1}^u\e^{-n(u-v)\sum_{k=1}^d|z_k|^{\alpha_k}}\d v\nonumber
 \\&&\qquad\leq\frac{Kn}{F_n^2}\int_{\R^d} \frac{|\widehat{\phi}(z)|^2}{\sum_{k=1}^d|z_k|^{\alpha_k}}\d
 z\int_{t_1}^{t_2}\e^{-n\delta_n u}\frac{1-\e^{-n(u-t_1)\sum_{k=1}^d|z_k|^{\alpha_k}}}{\sum_{k=1}^d|z_k|^{\alpha_k}}\d
 u.
 \eeqlb
Since $t_2\geq
 t_1\geq\varepsilon$, using the inequality $1-\e^{-x}\leq x^r$ for all $x\geq0$ and $r\in(0, 1]$, from (\ref{s3-5-2}) we get that
  \beqnn
 I_{11}(n,\psi_n)&\leq&\frac{Kn^{1+r}}{F_n^2}\e^{-n\delta_n\varepsilon}
 \int_{\R^d} \frac{|\widehat{\phi}(z)|^2}{(\sum_{k=1}^d|z_k|^{\alpha_k})^{2-r}}\d
 z\int_{t_1}^{t_2}(u-t_1)^{r}\d
 u,
 \eeqnn
for every $r\in(0, 1)$. Lemma
\ref{s2-lm-1} implies $\int_{\R^d}
\frac{|\widehat{\phi}(z)|^2}{(\sum_{k=1}^d|z_k|^{\alpha_k})^{2-r}}\d
 z<\infty$. In addition,
 $n^{\kappa}\delta_n\to\theta\in(0,\infty)$ with $\kappa\in(0, 1)$ implies that for all
 $n$, $\frac{Kn^{1+r}}{F_n^2}\e^{-n\delta_n\varepsilon}$ are bounded. Consequently, there
 exists a constant $K$ independent of $t_1$ and $t_2$, such that
  \beqlb\label{s3-5-4}
  I_{11}(n,\psi_n)&\leq&K|t_2-t_1|^{1+r}.
 \eeqlb

 {\bf (II)} We then proceed to estimate $\tilde{I}_2(n,\psi_n)$. Since $\bar{f}_n$ is bounded,
 applying (\ref{s3-6}), (\ref{s3-9})
and (\ref{s1-2}) to (\ref{s3-4-1})  we
obtain that for some constant $K>0$,
 \beqnn
 \tilde{I}_2(n,\psi_n)\leq\frac{K}{F_n^2}\int_0^n\e^{-\delta_n s}\tilde{h}(\frac{s}{n})\d
 s\int_0^{n-s}\e^{-\delta_n u}\tilde{h}(\frac{s+u}{n})\d
 u\int_{\R^d} \phi(x)T_u\phi(x)\d x.
 \eeqnn
Then by the same arguments as led to (\ref{s3-5-2}), we have that
 \beqlb\label{s3-5-5}
 \tilde{I}_2(n,\psi_n)&\leq&\frac{Kn^2}{F_n^2}\int_{t_1}^{t_2}\e^{-n\delta_n s}\d
 s\int_{s}^{t_2}\e^{-n\delta_n (u-s)}\d u
 \int_{\R^d} |\widehat{\phi}(z)|^2\e^{-n(u-s)\sum_{k=1}^d|z_k|^{\alpha_k}}\d z\nonumber
 \\&\leq&\frac{Kn}{F_n^2}\int_{\R^d} \frac{|\widehat{\phi}(z)|^2}{\sum_{k=1}^d|z_k|^{\alpha_k}}\d
 z\int_{t_1}^{t_2}\e^{-n\delta_n s}(1-\e^{-n(t_2-s)\sum_{k=1}^d|z_k|^{\alpha_k}})\d s.
 \eeqlb
Consequently, repeating the same arguments used to (\ref{s3-5-4}), we can readily get
that for any $r\in (0, 1)$, there
exists a constant $K>0$ independent of $t_1$ and $t_2$ such that
  \beqlb\label{s3-5-7}
  \tilde{I}_2(n,\psi_n)&\leq&K|t_2-t_1|^{1+r}.
 \eeqlb

 {\bf (III)} At last, we consider $\tilde{I}_3(n,\psi_n)$. Since $\delta_n\to0$, without loss of generality, we can assume $\delta_n<\gamma$.
  Using again (\ref{s1-2}), (\ref{s3-6}), (\ref{s3-9}), (\ref{s2-3})
and (\ref{s2-6}) to (\ref{s3-2-3}), from the condition (\ref{s3-5-1}), we
get that for any $r\in(0, 1)$,
  \beqlb\label{s3-5-8}
   \tilde{I}_3(n,\psi_n)&\leq&\frac{\delta_n}{F_n^2}\int_0^n\e^{-\delta_n s}\d s\int_0^{n-s}\e^{-\gamma u}\d u
  \int_0^u \tilde{h}(\frac{s+v}{n})\tilde{h}(\frac{s+u}{n})\d v\int_{\R^d}|\widehat{\phi}(z)|^2\d z\nonumber
  \\&=&\frac{n^3\delta_n}{F_n^2}\int_{\R^d}|\widehat{\phi}(z)|^2\d z\int_0^1\e^{-n\delta_n s}
  \d s\int_s^{1}\e^{-n\gamma(u-s)}\tilde{h}(u)\d u
  \int_s^u \tilde{h}(v)\d v\nonumber
  \\&=&\frac{n^3\delta_n}{F_n^2}\int_{\R^d}|\widehat{\phi}(z)|^2\d z\int_{t_1}^{t_2}\e^{-n\gamma u}\d u\int_{t_1}^u\d v
  \int_0^v\e^{ns(\gamma-\delta_n)}\d s\nonumber
  \\&\leq&\frac{n\delta_n}{(\gamma-\delta_n)^2F_n^2}\int_{\R^d}|\widehat{\phi}(z)|^2\d z
  \int_{t_1}^{t_2}\e^{-n\delta_n t_1}(1-\e^{-n\gamma(u-t_1)})\d
  u\nonumber
  \\&\leq&\frac{\gamma^rn^{1+r}\delta_n}{(\gamma-\delta_n)^2F_n^2}\e^{-n\delta_n \varepsilon}\int_{\R^d}|\widehat{\phi}(z)|^2\d
  z|t_2-t_1|^{1+r}.
 \eeqlb
By the same reason as applied to (\ref{s3-5-4}), for every $r\in(0,
 1)$ there exists a constant $K$ independent of $t_1$ and $t_2$ such that
 \beqlb\label{s3-5-9}
 \tilde{I}_3(n,\psi_n)\leq K|t_2-t_1|^{1+r}.
 \eeqlb

 Summing up, from (\ref{s3-5-4}),(\ref{s3-5-7}) and (\ref{s3-5-9}) we know that for any
$r\in(0, 1)$, there is a constant $K$ which is
independent of $t_1, t_2$ such that
 \beqlb\label{s4-1-3}
 |\tilde{I}_3(n, i\omega\psi_n)|+|\tilde{I}_2(n, i\omega\psi_n)|+|I_{11}(n, i\omega\psi_n)|\leq K(\phi, r)\omega^2|t_2-t_1|^{1+r}.
 \eeqlb
Note that
 \beqlb\label{s4-1-4}
 &&\Big|1-{\rm Re}\Big(\bfE\big(\exp\{-i\omega\langle\tilde{X}_n,
 \phi h\rangle\}\big)\Big)\Big|\nonumber
 \\&&\qquad\qquad\leq |I_1(n, i\omega\psi_n)|+|I_2(n, i\omega\psi_n)|+|I_3(n,
 i\omega\psi_n)|.\qquad
 \eeqlb
Therefore, (\ref{s4-1-2}), (\ref{s4-1-3}) and (\ref{s4-1-4}) yield that
 $$\int_0^{1/\eta}\Big(1-{\rm Re}\big(\bfE(\exp\{-i\omega\langle\tilde{X}_n,
 \phi h\rangle\})\big)\Big)\d\omega\leq\frac{K(\phi,
 r)}{3\eta^3}|t_2-t_1|^{1+r},$$
which completes the proof of (\ref{s4-1-1}) and hence the
proof of Theorem 2.1.\qed\medskip

\noindent{\bf Acknowledgments}

Thanks due to the anonymous referees
for careful reading of the paper and for their
suggestions which have helped to improve the quality of the paper.


\end{document}